\documentclass[12pt]{article}
\usepackage[utf8]{inputenc}
\usepackage[english]{babel}
\usepackage{amsmath}
\usepackage{amsfonts}
\usepackage{graphicx}
\usepackage{fouridx}
\usepackage{amsthm}
\usepackage{amssymb}
\usepackage{amscd}
\usepackage{enumerate}
\usepackage{mathtools}
\usepackage[french]{varioref}
\usepackage{multirow} %tableau
\usepackage{hyperref}
\usepackage{bm}
\usepackage{multicol}
\usepackage{multirow}
\usepackage{xcolor}

\usepackage{cite}

\usepackage{hyperref}
\hypersetup{
colorlinks   = true,
citecolor    = black,
linkcolor =black,
urlcolor=black,
}

\hypersetup{hidelinks}

\usepackage{enumitem}  % http://ctan.org/pkg/enumitem
\usepackage{calc}% http://ctan.org/pkg/calc
\setlist{labelindent=1pt,itemsep=.5em}
\setlist[itemize]{leftmargin=1.2cm}
\setlist[enumerate]{itemindent=0em,leftmargin=1.2cm}
\setlist[enumerate,1]{label={\upshape(\roman*)}}

\newtheorem{theorem}{Theorem}[section]

\newtheorem{thm}[theorem]{Theorem}
\newtheorem{lem}[theorem]{Lemma}
\newtheorem{prop}[theorem]{Proposition}
\newtheorem{example}{Example}[section]
\newtheorem{defn}[theorem]{Definition}
\newtheorem{remq}[theorem]{Remark}
\numberwithin{equation}{section}
\usepackage{multicol}
\input{epsf.sty}
\usepackage{varioref}

\newcommand{\N}{\mathbb{N}}
\newcommand{\C}{\mathbb{C}}

\newcommand{\K}{\mathbb{K}}

\allowdisplaybreaks

\makeatletter
\newcommand{\subjclass}[2][2010]{%
  \let\@oldtitle\@title%
  \gdef\@title{\@oldtitle\footnotetext{#1 \emph{Mathematics subject classification}: #2}}%
}
\newcommand{\keywords}[1]{%
  \let\@@oldtitle\@title%
  \gdef\@title{\@@oldtitle\footnotetext{\emph{Keywords}: #1}}%
}
\makeatother

\title{Classification, centroids and derivations of two-dimensional Hom-Leibniz algebras}


\date{\today}

\author{\sc Anja Arfa{$^{1,2}$}, Nejib Saadaoui{$^{3}$}, Sergei Silvestrov{$^{4}$}}
\date{{\small{$^{1}$ Department of Mathematics, College of Sciences and Arts in Gurayat, Jouf University, Sakakah, Saudi Arabia. \text{amarfa@ju.edu.sa} \\
\small{$^{2}$    Faculty of Sciences, University of Sfax, BP 1171, 3000 Sfax, Tunisia\\
$^{3}$ Universit\'{e} de Gab\`{e}s, Institut Sup\'{e}rieur d'Informatique et de Multim\'{e}dia de  Gab\`{e}s, Campus universitaire-BP 122, 6033 Cit\'{e} El Amel 4, Gab\`{e}s, Tunisie.
\text{najib-saadaoui@ yahoo.fr}  \\
\small{$^{4}$}
Division of Mathematics and Physics, School of Education, Culture and
Communication, M\"{a}lardalen University, Box 883, 72123 V\"{a}ster{\aa}s, Sweden. \text{sergei.silvestrov@mdh.se}
 }}}}

\begin{document}
\maketitle
\begin{abstract}
Several recent results concerning Hom-Leibniz algebra are reviewed, the notion of symmetric Hom-Leibniz superalgebra is introduced and some properties are obtained. Classification of $2$-dimensional Hom-Leibniz algebras is provided. Centroids and derivations of multiplicative Hom-Leibniz algebras are considered including the detailed study of $2$-dimensional Hom-Leibniz algebras.
	\end{abstract}

\section*{Introduction}	
Hom-Lie algebras and quasi-Hom-Lie algebras were introduced first in 2003 in \cite{HLSPrepr2003JA2006:deformLiealgsigmaderiv} where a general method for construction of deformations and discretizations of Lie algebras of vector fields based on twisted derivations obeying twisted Leibniz rule was developped motivated by the examples of $q$-deformed Jacobi identities in $q$-deformations of Witt and Visaroro and in related $q$-deformed algebras discovered in 1990'th in string theory, vertex models of conformal field theory, quantum field theory and quantum mechanics, and also in development of $q$-deformed differential calculi and $q$-deformed homological algebra
\cite{
aizawasato199091:qdefViralgcenext,
ChaiElinPop1990:qconfalgcentext,
ChaiKuLuk,
ChaiPopPres,
ChaiIsLukPopPresn1991:Viralgconfdim,
daskaloyannis1992generalized,
Hu1992:qWittalgqLie,
Kassel92,
LiuKQuantumCentExt,
LiuKQ1992:CharQuantWittAlg,
LiuKQPhDthesis}.
The central extensions and cocycle conditions for general quasi-Hom-Lie algebras and Hom-Lie algebras, generalizing in particular $q$-deformed Witt and Virasoro algebras, have been first considered in \cite{HLSPrepr2003JA2006:deformLiealgsigmaderiv,LarssonSilvJA2005:QuasiHomLieCentExt2cocyid} and for graded color quasi-Hom-Lie algebras in \cite{SigurdssonSilvestrov2009:colorHomLiealgebrascentralext}.
General quasi-Lie and quasi-Leibniz algebras introduced in \cite{LarssonSilvestrov2005:QuasiLiealgebras} and color quasi-Lie and color quasi-Leibniz algebras introduced in \cite{LarssonSilvestrov2005:GradedquasiLiealgebras} in 2005,
include the Hom-Lie algebras, the quasi-Hom-Lie algebras as well as the color Hom-Lie algebras, quasi-Hom-Lie color algebras, quasi-Hom-Lie superalgebras and Hom-Lie superalgebras, and color quasi-Leibniz algebras, quasi-Leibniz superalgebras, quasi-Hom-Leibniz superalgebras and Hom-Leibniz algebras.

Hom-Lie algebras and Hom-Lie superalgebras and more general
color quasi-Lie algebras interpolate on the fundamental level of defining identities between Lie algebras, Lie superalgebras, color Lie algebras and related non-associative structures and their deformations, quantum deformations and discritizations, and thus might be useful tool for unification of methods and models of classical and quantum physics, symmetry analysis and non-commutative geometry and computational methods and algorithms based on general non-uniform discretisations of differential and integral calculi.
Binary Hom-algebra structures typically involve a bilinear binary operation and one or several linear unary operations twisting the defining identities of the structure in some special nontrivial ways, so that the original untwisted algebraic structures are recovered for the specific twisting linear maps. In quasi-Lie algebras and quasi-Hom-Lie algebras, the Jacobi identity contains in general six triple bracket terms twisted in special ways by families of linear maps, and the skew-symmetry is also in general twisted by a family of linear maps. Hom-Lie algebras is a subclass of quasi-Lie algebras with the bilinear product satisfying the Jacobi identity containing only three triple bracket terms twisted by a single linear map and the usual non-twisted skew-symmetry identity. Hom-Leibniz algebras arise when skew-symmetry is not required, while the Hom-Jacobi identity is written as Hom-Leibniz identity. Lie algebras and Leibniz algebras as a special case of Hom-Leibniz algebras are obtained for the trivial choice of the twisting linear map as the identity map on the underlying linear space. The

Hom-Lie algebras, Hom-Lie superalgebras, Hom-Leibniz algebras and Hom-Leibniz superalgebras with twisting linear map different from the identity map, are rich and complicated algebraic structures with classifications, deformations, representations, morphisms, derivations and homological structures being fundamentally dependent on joint properties of the twisting maps as a unary operations and bilinear binary product intrinsically linked by Hom-Jacobi or Hom-Leibniz identities.
Hom-Lie admissible algebras have been considered first in \cite{MakhoufSilvestrov:Prep2006JGLTA2008:homstructure}, where
the Hom-associative algebras and more general $G$-Hom-associative algebras including the Hom-Vinberg algebras (Hom-left symmetric algebras), Hom-pre-Lie algebras (Hom-right symmetric algebras), and some other new Hom-algebra structures have been introduced and shown to be Hom-Lie admissible, in the sense that the operation of commutator as new product in these Hom-algebras structures yields Hom-Lie algebras. Furthermore, in \cite{MakhoufSilvestrov:Prep2006JGLTA2008:homstructure},
flexible Hom-algebras and Hom-algebra generalizations of derivations and of adjoint derivations maps have been introduced, and the Hom-Leibniz algebras appeared for the first time, as an important special subclass of quasi-Leibniz algebras introduced in \cite{LarssonSilvestrov2005:QuasiLiealgebras} in connection to general quasi-Lie algebras following the standard Loday’s conventions for Leibniz algebras (i.e. right Loday algebras) \cite{LodayEM1993:versnoncomalgLie}.
In \cite{MakhoufSilvestrov:Prep2006JGLTA2008:homstructure}, moreover the investigation of classification of finite-dimensional Hom-Lie algebras have been initiated with construction of families of the low-dimensional Hom-Lie algebras.
Since \cite{HLSPrepr2003JA2006:deformLiealgsigmaderiv,
LarssonSilvJA2005:QuasiHomLieCentExt2cocyid,
LarssonSilvestrov2005:GradedquasiLiealgebras,
LarssonSilvestrov2005:QuasiLiealgebras,
LarssonSilvestrov:Quasidefsl2,
MakhoufSilvestrov:Prep2006JGLTA2008:homstructure}, Hom-algebra structures expanded into a popular area providing a new broad framework for establishing fundamental links between deformations and quantum deformations of associative algebras, various classes of non-associative algebras, super-algebras, color algebras, $n$-ary algebraic structures, non-commutative differential calculus and homological algebra constructions for associative and non-associative structures.
Quadratic Hom-Lie algebras have been considered in \cite{benayadimakhloufJGP2014:homliesyminvnondegbilform} and
representation theory, cohomology and homology theory of Hom-Lie algebras have been considered in
\cite{AmmarEjbehiMakhlouf2011:cohomhomdeformhomalg,
AmmarMabroukMakhloufCohomnaryHNLalg2011,sheng2012representations,Yau2009:HomolHom}.
Investigation of color Hom-Lie algebras and Hom-Lie superalgebras and $n$-ary generalizations expanded
\cite{
AbAmMakh:HomaltHomMalcHomJordSuperalg,
AbdaouiAmmarMakhloufCohhomLiecolalg2015,
AbramovSilvestrov:3homLiealgsigmaderivINvol,
AmmarMakhloufJA2010:homliesuperaladmsuperalg,
AmmarMakhloufSaadaoui2013:CohlgHomLiesupqdefWittSup,
AmmarMakhloufSilv:TernaryqVirasoroHomNambuLie,
AmAyMabMakh:QuadrColHomLieAlgs,
ArmakanFarhangdoost:2017IJGMMP:GeomaspectsextHomLiesuperalgs,
ArmakanFarhSilv2021:ndKilformsHomLiesuper,
armakanrazavicomalg2020:completehomliesuper,
ArmakanSilv2020:envelalgcertaintypescolorHomLie,
ArmakanSilvFarh:envelopalgcolhomLiealg,
ArmakanSilvFarh20172019:exthomLiecoloralg,
ArmakanSilv:colorHomLieHomLiebnOmniHomLie,
akms:ternary,
ams:ternary,
ArnlindMakhloufSilvnaryHomLieNambuJMP2011,
AtMaSi:GenNambuAlg,
BeitesKaygorodovPopov2018:GenDermultnaryHomOmegacolalg,
Bakayoko2014:ModulescolorHomPoisson,
BakayokoDialo2015:genHomalgebrastr,
BakyokoSilvestrov:Homleftsymmetriccolordialgebras,
BakyokoSilvestrov:MultiplicnHomLiecoloralg,
BakayokoToure2019:genHomalgebrastr,
CaoChen2012:SplitregularhomLiecoloralg,
GuanChenSun:HomLieSuperalgebras,
KitouniMakhloufSilvestrov:nhominduced,
kms:narygenBiHomLieBiHomassalgebras2020,
kms:solvnilpnhomlie2020,
liu2013hom,
MabroukNcibSilvestrov2020:GenDerRotaBaxterOpsnaryHomNambuSuperalgs,
Makhlouf2010:paradigmnonasshomalghomsuper,
MandalMishra:HomGerstenhaberHomLiealgebroids,
MishraSilvestrov:SpringerAAS2020HomGerstenhalgsHomLiealgds,
SigSilv:CzechJP2006:GradedquasiLiealgWitt,
SigurdssonSilvestrov2009:colorHomLiealgebrascentralext,
SilvestrovParadigmQLieQhomLie2007,
WangZhangWei2015:HomLeibnizsuperalg,
Yuan2012:HomLiecoloralgstr,
ZhouChenMa:GenDerHomLiesuper}.

At the same time, in recent years, the theory and classification of Leibniz algebras, and also  Leibniz superalgebras extending Leibniz algebras in a similar way as Lie superalgebras generalize Lie algebras, continued being actively investigated motivated in part by applications in Physics and by graded homological algebra structures of non-commutative differential calculi and non-commutative geometry \cite{Camacho,Hu}. In \cite{SaidHusainRakhimovBasri2017:centrdsderlowdimleibnizalgs}, the description of centroid and derivations of low-dimensional Leibniz algebras using classification results is introduced.
Note that a classification of two dimensional Leibniz algebras have been given by Loday in \cite{LodayEM1993:versnoncomalgLie}.
In dimension three there are fourteen isomorphism classes, and the list can be found in \cite{CasasILL:laa2012clas3dleibn,RikhsiboevRakhimovIJMP2012:clas3dLeibniz}. Furthermore, a classification of low-dimensional complex solvable Leibniz algebras can be found in \cite{CaneteKhudoyberdiyevLAA2013:clas4dLeibniz,
CasasILL:laa2012clas3dleibn,
CasasLOK:lma2013clasolvLeibniznullfilinilrad,
KhudLadraOmirov2014:solvLeibniznilrad,
KhudRakhimovSaidHusain2014:clas5dsolvLeibniz,
KhudLadraOmirov2014:solvLeibniznilrad}, and of two and three-dimensional complex Leibniz algebra is given in \cite{LodayEM1993:versnoncomalgLie,
BenkartNeherJA2006:CantroidExtAffinerootgrad}.

The purpose of the present work is to study the classification of multiplicative $2$-dimen\-sional Hom-Leibniz algebras and to investigate centroids and derivations of Hom-Leibniz algebras and superalgebras and the concepts of left, right and symmetric (two-sided) Hom-Leibniz superalgebras.
In Section \ref{sec:homLeibnizsuperalgs}, the left, right and symmetric (two-sided) Hom-Leibniz algebras and superalgebras generalizing the well known left, right and symmetric (two-sided) Leibniz algebras are defined and some of their properties are reviewed.
In Section \ref{sec:clasmult2dimhomLeibnizalgs}, we provide classification of multiplicative
$2$-dimensional Hom-Leibniz algebras.
In Section \ref{sec:centroidsderivhomLeibnizsuperalgs}, we review centroids and derivations of  Hom-Leibniz superalgebras and some of their properties.
In Section \ref{sec:centroidsderiv2dimmulthomLeibnizalgs},
we describe the algorithm to find centroids and derivations of Hom-Leibniz algebras,
apply it to the classification of $2$-dimensional Hom-Leibniz algebras from the previous section and
use properties of the centroids of Hom-Leibniz algebras to categorize the algebra into having small and not small centroids.

\section{Hom-Leibniz algebras and superalgebras}
\label{sec:homLeibnizsuperalgs}
%%%%%%%%%%%%%%%%%%%%%%%%%%%%%%%%%%%%%%
For exposition clarity, we assume throughout this article that all linear spaces are over a field $\mathbb{K}$ of characteristic different from $2$, and just note that many of the results in this article hold as formulated or just with minor modifications for any field. Multilinear maps $f\colon V_1 \times \cdots \times V_n \to W $ on finite direct products and linear maps $ F\colon V_1 \otimes \cdots \otimes V_n \to W\text{,}$ on finite tensor products of linear spaces are identified via $F(v_1\otimes \cdots \otimes v_n) = f(v_1,\ldots,v_n)$.
A linear space $V$ is said to be a $G$-graded by an abelian group $G$ if $V=\bigoplus\limits_{g\in G} V_g$ for a family $\{V_g\}_{g\in G}$ of linear subspaces of $V$.
For each $g\in G$, the elements of $V_g$ are said to be homogeneous of degree $g\in G$, and the set of homogeneous elements is the union $\mathcal{H}(V)= \bigcup\limits_{g\in G} V_g$ all the spaces $V_g$ of homogenous elements of degree $g$ for all $g\in G$.
For two $G$-graded linear spaces $V=\bigoplus\limits_{g\in G} V_g$ and $V'=\bigoplus\limits_{g\in G} V'_g$, a linear mapping $f : V\rightarrow V'$ is called homogeneous of degree $d$ if $f(V_g)\subseteq  V'_{g+d}$, for all $g\in G.$
The homogeneous linear maps of degree zero (even maps) are those homogeneous linear maps satisfying  $f(V_g)\subseteq V'_{g}$ for any $g\in G$.
In the $\mathbb{Z}_{2}$-graded linear spaces $V=V_{0}\oplus V_{1}$, also called superspaces,
the elements of $V_{|j|}$, $|j|\in \mathbb{Z}_{2},$ are said to be homogenous of parity $|j|.$
The space $End (V)$ is $\mathbb{Z}_{2}$-graded with a direct sum $End (V)=(End (V))_{0}\oplus(End (V))_{1}$ where
$(End (V))_{|j|}=\{f\in End (V) \mid f (V_{i})\subset V_{i+j}\}.$
The elements of $(End (V))_{}$  are said to be homogenous of parity $|j|.$

An algebra $(A, \cdot)$ is called $G$-graded if its linear space is $G$-graded $A=\bigoplus\limits_{g\in G}A_g$, and $A_g\cdot A_h\subseteq A_{g+h}$ for all $g, h\in G$. Homomorphisms of $G$-graded algebras $A$ and $A'$ are homogeneous algebra morphisms of degree $0_G$ (even maps).
Hom-superalgebras are triples $(A, \mu, \alpha)$ consisting of a $\mathbb{Z}_2$-graded linear space (superspace) $A=A_0\oplus A_1$, an even bilinear map $\mu : A\times A\rightarrow A$, and
an even linear map $\alpha : A\rightarrow A$.
An even linear map $f : A \rightarrow A'$ is said to be a weak morphism of hom-superalgebras if it is algebra structures homomorphism
($f\circ\mu=\mu'\circ(f\otimes f)$), and a morphism of hom-superalgebras if moreover  $f\circ\alpha=\alpha'\circ f$. The important point for understanding of difference between algebras and Hom-algebras is that the properties preserved by the morphisms of specific Hom-algebras do not need to be preserved by all weak morphisms between these Hom-algebras, and the classifications of Hom-algebras up to weak isomorphisms (all algebra isomorphisms) and of Hom-algebras up to isomorphisms of Hom-algebras differ substantially in that the set of isomorphisms intertwining the twisting maps $\alpha$ and $\alpha'$ often is a proper subset of all isomorphisms, and thus classification of Hom-algebras up to isomorphism of Hom-algebras of different types typically contains much more isomorphism classes than weak isomorphism classes.

In any Hom-superalgebra $(A=A_0\oplus A_1, \mu, \alpha)$,
\begin{align*}
\mu(A_0, A_0)\subseteq A_0, \quad \mu(A_1, A_0)\subseteq A_1, \quad
\mu(A_0, A_1)\subseteq A_1, \quad \mu(A_1, A_1)\subseteq A_0.
\end{align*}
Hom-subalgebras (graded Hom-subalgebras) of Hom-superalgebra $(A,\mu,\alpha)$ are   $\mathbb{Z}_{2}$-graded linear subspaces $I =  (I\cap A_{0}) \oplus (I\cap A_{1})$ of $A$ obeying   $\alpha(I) \subseteq I$ and $\mu(I,I) \subseteq I$.
Hom-associator on Hom-superalgebra $(A=A_0\oplus A_1, \mu, \alpha)$ is the even trilinear map   $as_{\alpha,\mu}=\mu\circ(\mu\otimes\alpha-\alpha\otimes\mu):A\times A\times A\rightarrow A,$
acting on elements as
$as_{\alpha,\mu}(x, y, z)=\mu(\mu(x, y), \alpha(z))-\mu(\alpha(x), \mu(y, z)),$
or $as_{\alpha,\mu}(x, y, z)=(xy)\alpha(z)-\alpha(x)(yz)$ in juxtaposition notation $xy=\mu(x, y)$.
Since $|as_{\alpha,\mu}(x,y,z))| =  |x|+|y|+|z|$ for $x, y, z \in \mathcal{H}(A)=A_0\cup A_1$ in any Hom-superalgebra $(A=A_0\oplus A_1,\mu, \alpha),$
\begin{align}
as_{\alpha,\mu}(A_0, A_0, A_0)\subseteq A_0, & \quad
as_{\alpha,\mu}(A_1, A_1, A_0)\subseteq A_0, \\
as_{\alpha,\mu}(A_1, A_0, A_1)\subseteq A_0, & \quad
as_{\alpha,\mu}(A_0, A_1, A_1)\subseteq A_0, \\
as_{\alpha,\mu}(A_1, A_0, A_0)\subseteq A_1, & \quad
as_{\alpha,\mu}(A_0, A_1, A_0)\subseteq A_1, \\
as_{\alpha,\mu}(A_0, A_0, A_1)\subseteq A_1, & \quad
as_{\alpha,\mu}(A_1, A_1, A_1)\subseteq A_1.
\end{align}

\begin{defn}[\cite{HLSPrepr2003JA2006:deformLiealgsigmaderiv,MakhoufSilvestrov:Prep2006JGLTA2008:homstructure}] \label{def:homlie}
Hom-Lie algebras  are triples $(\mathcal{G}, [\cdot,\cdot], \alpha)$ consisting of a linear space $\mathcal{G}$ over a field $ \mathbb{K}$, a  bilinear map $[\cdot,\cdot]:\mathcal{G}\times \mathcal{G}\rightarrow \mathcal{G}$ and a $ \mathbb{K}$-linear map
$ \alpha:\mathcal{G}\rightarrow \mathcal{G} $
satisfying for all $x,\, y,\, z \in \mathcal{G}$,
\begin{align}
[x, y]&=-[y,x], & \textup{Skew-symmetry} \label{skewsym} \\
[\alpha (x), [y, z]]+ [\alpha (y), [z, x]] &+ [\alpha (z), [x, y]] = 0. & \textup{Hom–Jacobi identity}
\label{homliejacobi}
	\end{align}
\end{defn}

Hom-Lie algebra is called a multiplicative Hom-Lie algebra if $\alpha$ is an algebra morphism, $\alpha([\cdot,\cdot]) = ([\alpha(.),\alpha(.)])$, meaning that $\alpha([x,y])=[\alpha(x),\alpha(y)]$ for any $x, y \in \mathcal{G}$. Lie algebras are a very special subclass of multiplicative Hom-Lie algebras obtained for $\alpha = id$ in Definition \ref{def:homlie}.

\begin{defn}[\cite{AmmarMakhloufJA2010:homliesuperaladmsuperalg,LarssonSilvestrov2005:GradedquasiLiealgebras}] \label{def:homliesuper} Hom-Lie superalgebras are triples $(\mathcal{G},[\cdot,\cdot],\alpha)$ which consist of $\mathbb{Z}_{2}$-graded linear space $\mathcal{G} = \mathcal{G}_{0} \oplus \mathcal{G}_{1}$, an even bilinear map  $[\cdot,\cdot] : \mathcal{G} \times \mathcal{G} \to \mathcal{G}$ and an even linear map $\alpha : \mathcal{G} \to \mathcal{G}$ satisfying the super skew-symmetry and Hom-Lie super Jacobi identities for homogeneous elements $x, y, z \in \mathcal{H}(\mathcal{G})$,
\begin{align}
[x,y] =-(-1)^{|x||y|}[y,x], \quad \quad \quad  \textup{Super skew-symmetry}   \label{superskewsymmetry} \\
(-1)^{|x||z|}[\alpha(x),[y,z]] +(-1)^{|y||x|}[\alpha(y),[z,x]] +(-1)^{|z||y|}[\alpha(z),[x,y]]=0. \label{homliesuperjacobi}  \\
\textup{Super Hom-Jacobi identity} \nonumber
\end{align}
\end{defn}
Hom-Lie superalgebra is called multiplicative Hom-Lie superalgebra, if $\alpha$ is an algebra morphism, $\alpha([x,y])=[\alpha(x),\alpha(y)]$ for any $x, y \in \mathcal{G}$.

\begin{remq}
In any Hom-Lie superalgebra, $(\mathcal{G}_{0},[\cdot,\cdot],\alpha)$ is a Hom-Lie algebra since
$[\mathcal{G}_{0},\mathcal{G}_{0}]\in \mathcal{G}_{0}$ and $\alpha(\mathcal{G}_{0})\in \mathcal{G}_{0}$ and $(-1)^{|a||b|}= (-1)^{0}=1$ for $a,b \in \mathcal{G}_{0}$. Thus,
Hom-Lie algebras can be also seen as special class of Hom-Lie superalgebras when $\mathcal{G}_{1}=\{0\}$.
\end{remq}
\begin{remq}
From the point of view of Hom-superalgebras, Lie superalgebras is a special subclass of multiplicative Hom-Lie superalgebras obtained when $\alpha = id$ in Definition \textup{\ref{def:homliesuper}} which becomes the definition of Lie superalgebras as  $\mathbb{Z}_{2}$-graded linear spaces $\mathcal{G} = \mathcal{G}_{0} \oplus \mathcal{G}_{1}$, with a graded Lie bracket $[\cdot,\cdot] : \mathcal{G} \times \mathcal{G} \to \mathcal{G}$ of degree zero, that is $[\cdot,\cdot]$ is a bilinear map obeying
$[\mathcal{G}_{i}, \mathcal{G}_{j}] \subset \mathcal{G}_{{i+j}(mod2)},$
and for $x, y, z \in \mathcal{H}(\mathcal{G})=\mathcal{G}_{0} \cup \mathcal{G}_{1}$,
\begin{align}
[x,y]=-(-1)^{|x||y|}[y,x],  \quad \quad \quad \textup{Super skew-symmetry}   \\
(-1)^{|x||z|}[x,[y,z]] +(-1)^{|y||x|}[y,[z,x]] +(-1)^{|z||y|}[z,[x,y]]=0. \label{liesuperjacobi}
 \begin{array}[c]{c}
\textup{Super Jacobi} \\
\textup{identity}
\end{array}
\end{align}
In super skew-symmetric superalgebras, the super Hom-Jacobi identity can be presented equivalently in the form of super Leibniz rule for the maps $ad_x=[x,\cdot]:\mathcal{G} \rightarrow \mathcal{G},$
\begin{equation}
[x,[y,z]]=[[x,y],z]+(-1)^{|x||y|}[y,[x,z]]. \label{liesuperjacobi:leibniz}
\end{equation}
Hom-Lie superalgebras are substantially different from Lie superalgebras, as all algebraic structure properties, morphisms, classifications and deformations become  fundamentally dependent on the joint structure and properties of unary operation given by the linear map $\alpha$ and the bilinear product $[\cdot,\cdot]$ intricately linked via the $\alpha$-twisted super-Jacoby identity \eqref{homliesuperjacobi}.
\end{remq}

In super skew-symmetric Hom-superalgebras, the super Hom-Jacobi identity can be presented equivalently in the form of super Hom-Leibniz rule for the maps $ad_x=[x,.]:\mathcal{G} \rightarrow \mathcal{G},$
\begin{equation}
[\alpha(x),[y,z]]=[[x,y],\alpha(z)]+(-1)^{|x||y|}[\alpha(y),[x,z]] \label{homliesuperjacobi:leibniz}
\end{equation}
since, by super skew-symmetry \eqref{superskewsymmetry}, the following equalities are equivalent:
\begin{align*}
& [\alpha(x),[y,z]] = [[x,y],\alpha(z)]+(-1)^{|x||y|}[\alpha(y),[x,z]], \\
& [\alpha(x),[y,z]]-[[x,y],\alpha(z)]-(-1)^{|x||y|}[\alpha(y),[x,z]]=0, \\
& [\alpha(x),[y,z]]+(-1)^{|z|(|x|+|y|)}[\alpha(z),[x,y]]-(-1)^{|x||y|}[\alpha(y),[x,z]]=0, \\
& [\alpha(x),[y,z]]+(-1)^{|z|(|x|+|y|)}[\alpha(z),[x,y]]-(-1)^{|x||y|}[\alpha(y),[x,z]]=0, \\
& [\alpha(x),[y,z]]+(-1)^{|z|(|x|+|y|)}[\alpha(z),[x,y]]+(-1)^{|x||y|}(-1)^{|z||x|}[\alpha(y),[z,x]]=0, \\
& (-1)^{|z||x|}[\alpha(x),[y,z]]+ (-1)^{|z||x|}(-1)^{|z|(|x|+|y|)}[\alpha(z),[x,y]]+(-1)^{|x||y|}[\alpha(y),[z,x]]=0, \\
& (-1)^{|x||z|}[\alpha(x),[y,z]]+(-1)^{|z||y|}[\alpha(z),[x,y]]+(-1)^{|y||x|}[\alpha(y),[z,x]]=0, \\
& (-1)^{|x||z|}[\alpha(x),[y,z]]+(-1)^{|z||y|}[\alpha(z),[x,y]]+(-1)^{|y||x|}[\alpha(y),[z,x]]=0.
\end{align*}

\begin{remq}
If skew-symmetry \eqref{skewsym} does not hold, then \eqref{homliesuperjacobi} and
\eqref{homliesuperjacobi:leibniz} are not necessarily equivalent, defining different Hom-superalgebra structures.
The Hom-superalgebras defined by just super algebras identity \eqref{homliesuperjacobi:leibniz} without requiring super Hom-skew-symmetry on homogeneous elements are Leibniz Hom-superalgebras, a special class of general $\Gamma$-graded quasi-Leibniz algebras (color quasi-Leibniz algebras) first introduced in \textup{\cite{LarssonSilvestrov2005:QuasiLiealgebras, LarssonSilvestrov2005:GradedquasiLiealgebras}}.
\end{remq}

%%%%%%%%%%%%%%%%%%%%%%%
%%%%%%%%%%%%%%%%%%%%%%%
\section{Symmetric (two-sided) Hom-Leibniz superalgebras}
\label{subsec:symetrichomLeibnizalgs}
In this section, we define the left, right and symmetric (two-sided) Hom-Leibniz algebras
Hom-Leibniz superalgebras and give examples of the symmetric Hom-Leibniz algebras and superalgebras.
We also study some properties of centroids of Hom-Leibniz superalgebras.

With the skew-symmetry \eqref{skewsym} satisfied, the Hom-Jacobi identity \eqref{homliejacobi} can be presented in two equivalent ways, for all $  x,\ y,\ z\in \mathcal{G}$,
\begin{alignat}{3}
\left[\alpha(x),[y,z]\right] &=	
\left[ [x,y],\alpha(z)\right] +  \left[\alpha(y), [x,z]\right], \quad && \text{(Left Hom-Leibniz)} \label{leftLeibniz}\\
	\left[\alpha(x),[y,z]\right] &=	
\left[ [x,y],\alpha(z)\right] - \left[ [x,z],\alpha(y)\right]. \quad && \text{(Right Hom-Leibniz)} \label{rigtLeibniz}
\end{alignat}
Without skew-symmetry however, these identities lead to different types of Hom-algebra structures both containing Hom-Lie algebras as subclass obeying the skew-symmetry \eqref{skewsym}.

\begin{defn}[\cite{MakhoufSilvestrov:Prep2006JGLTA2008:homstructure,LarssonSilvestrov2005:QuasiLiealgebras}]
Left Hom-Leibniz algebras are triples $(L, [\cdot,\cdot], \alpha)$ consisting of a linear space $L$ over a field $\K$, a bilinear map $[\cdot,\cdot]: L\times L\rightarrow L$ and
a linear map $\alpha:L\rightarrow L$ satisfying \eqref{leftLeibniz}
for all $  x,y,z\in L$.
Right Hom-Leibniz algebras are triples $(L, [\cdot,\cdot], \alpha)$ consisting of a linear space $L$ and over a field $\K$, a bilinear map $[\cdot,\cdot]: L\times L\rightarrow L$ and
a linear map $\alpha:L\rightarrow L$ satisfying \eqref{rigtLeibniz}
for all $  x, y, z\in L$.
\end{defn}

\begin{remq}
If $\alpha=id_L$, then  a left (right) Hom-Leibniz algebra is just a left (right) Leibniz algebra \textup{\cite{BenayadiHidri2016:LeibnizalgsinvbilinformsrelLiealgs,MasonaYamskulnabSigma2013:LeibnizLie}}.
Any Hom-Lie algebra is both a left Hom-Leibniz algebras and a right Hom-Leibniz algebra.
A left Hom-Leibniz algebra or a right Hom-Leibniz algebra is a Hom-Lie algebra if
$[x,x]=0$ for all $x\in L$. For the field $\K$ of characteristic different from $2$,
left Hom-Leibniz algebras and right Hom-Leibniz algebras are Hom-Lie algebras are Hom-Lie algebras if and only if $[x,x]=0$ for all $x\in L$.
\end{remq}

\begin{defn}\label{1thsymmetic }
A triple $(L,[\cdot,\cdot],\alpha)$ is called a symmetric (or two-sided) Hom-Leibniz algebra if it is both a left Hom-Leibniz algebra and a right Hom-Leibniz algebra, that is if both \eqref{leftLeibniz} and \eqref{rigtLeibniz} are satisfied for all $  x,\ y,\ z\in L$.
\end{defn}

\begin{remq}
	A left Hom-Leibniz algebra $(L,[\cdot,\cdot],\alpha,\beta)$  is a symmetric
	Hom-Leibniz algebra if and only if, for all  $x,\, y,\, z \in L$,
	\begin{equation}\label{thsymmetic}
	[ \alpha(y),[  x ,z  ] ] =-   [[  x ,z  ] ,\alpha (y ) ].
	\end{equation} 	
\end{remq}

\begin{example}
	Let $ (x_1,x_2,x_3) $ be a basis of $3$-dimensional space $  \mathcal{G}$ over $ \K $.
	% The following bracket and linear map $\alpha$ on $ \mathcal{G} $.
	% define a Hom-Lie algebra over $\K$:
	%	\begin{align*}
	%	[x_1,x_2]&=0    &\alpha(x_1)=a x_1+b x_2\\
	%	[x_1,x_3]&=x_1     &\alpha(x_2)=d x_1+e x_2\\
	%	[x_2,x_3]&= x_2    &\alpha(x_3)= x_3.\\
	%	\end{align*}
	Define a bilinear bracket operation on $\mathcal{G}\otimes \mathcal{G}  $ by
	%\[[x\otimes y,a\otimes b]=[[x,y],a]\otimes b+a\otimes[[x,y],b]\]
	\begin{align*}
	[x_1\otimes x_3,x_1\otimes x_3]&=  x_1\otimes x_1,\\	
[x_2\otimes x_3,x_1\otimes x_3]&=  x_2\otimes x_1,\\
	[x_2\otimes x_3,x_2\otimes x_3]&=  x_2\otimes x_2,
	\end{align*}
	with the other necessary brackets being equal to $ 0 $. For any linear map $\alpha$ on $ \mathcal{G} $,
	the triple $\left( \mathcal{G}\otimes \mathcal{G},[\cdot,\cdot],\alpha\otimes\alpha\right)$ is not
	a Hom-Lie algebra but it is a symmetric Hom-Leibniz algebra.
\end{example}
In the following examples, we construct Hom-Leibniz algebras on a linear space $ L\otimes L$ starting from a Lie  or a Hom-Lie algebra $L$.
%\begin{example}
%	The Jackson  $ \mathfrak{s}l(2) $ (see \cite{Ayedi} ) is defined with respect to a basis $ (x_1,x_2,x_3)$  by the brackets and a linear map $ \alpha $ such that 	
%	\begin{align*}
	%[x_1, x_2]&=  -2q x_{2}\qquad   &\alpha(x_1)&=qx_1\\
%	[x_1, x_3]&=2 x_3\qquad   &\alpha(x_2)&=q^2x_2\\
%	[x_2, x_3]&=-\frac{1}{2}(1+q)  x_1\qquad   &\alpha(x_3)&=qx_3
%	\end{align*}
%\end{example}

%In the sequel, Hom-Leibniz algebra refers to left or right or symmetric Hom-Leibniz algebra.
%Let $ (e_1,\cdots,e_n) $	a basis of linear space $ L $,
%$ [\cdot,\cdot]\, :L\times  L\rightarrow  L  $ a  bilinear map and  $ \alpha:L\rightarrow L \ $ a $ \mathbb{K}$-linear map. And set
%\begin{list}{$\bullet$}{}
%%%%%	\item $\displaystyle[e_i,e_j]=\sum_{r=1}^{n}B_r(e_i,e_j)e_r$, $B_r=\left(B_r(e_i,e_j)\right)_{\substack{1\leqslant  i\leqslant n\\1\leqslant j \leqslant n}}   $
%	\item $ p=(i,j)  $ with $ p=(i-1)n+k $ (Euclidean division).
%\item The matrix $ X=(x_{p,q}))_{\substack{1\leqslant  p\leqslant n\\1\leqslant q \leqslant n^2}}$\quad  where $q=(j-1)n+k$\quad	and \\$ x_{pq}=B_p(e_j,e_k)= x_{p,(j,k)} $
%	\item $\displaystyle\left[ \alpha(e_i),[e_j,e_k]\right] =\sum_{r=1}^{n}g_r(e_i,e_j,e_k)e_r$
%\end{list}

%%%%%%%%%%%%%%%%
\begin{prop}[\cite{KurdianiPirashviliJLT2002:Leibnizalgsstrsectensorpower}]
\label{tensorLie}
For any Lie algebra $(\mathcal{G}, [\cdot,\cdot])$, the bracket	 \[[x\otimes y,a\otimes b]=[x,[a,b]]\otimes y+x\otimes[y,[a,b]]\] defines a
 Leibniz algebra structure on the linear space $\mathcal{G}\otimes\mathcal{G}$.
\end{prop}
\begin{prop}[\cite{Yau:HomEnv}] \label{tensorLeibn}
	Let $(L,[\cdot,\cdot] )$ be a Leibniz algebra and $ \alpha : L\to  L$  be a Leibniz algebra endomorphism. Then  $(L,[\cdot,\cdot]_{\alpha},\alpha )$ is a Hom-Leibniz algebra, where $[x,y]_{\alpha}=[\alpha(x),\alpha(y)]  $.
\end{prop}
%\begin{remq}
%it's similar to the Hom-Lie algebra case given in \cite{Yau}.
%\end{remq}
Using Proposition \ref{tensorLeibn} and  Proposition \ref{tensorLie}, we obtain the following result.
\begin{prop}
	Let $(\mathcal{G},[\cdot,\cdot]') $ be a Lie algebra and $ \alpha : \mathcal{G}\to  \mathcal{G}$  be a Lie algebra endomorphism. We define on $\mathcal{G}\otimes\mathcal{G}$ the following bracket on $\mathcal{G}\otimes \mathcal{G},$
	\[ \left[ x\otimes y,a\otimes b\right] =\left[ \alpha(x),[\alpha(a),\alpha(b)]']'\otimes \alpha(y)+\alpha(x)\otimes[\alpha(y),[\alpha(a),\alpha(b)]'\right]'.  \]
Then
	$ \left(\mathcal{G}\otimes \mathcal{G},[\cdot,\cdot],\alpha \otimes \alpha \right)  $
	is a right  Hom-Leibniz algebra.
\end{prop}

\begin{prop}
	If  $(L,[\cdot,\cdot])$  is a symmetric Leibniz algebra,
	$\alpha \colon L\to  L$ a  morphism of Leibniz algebra, and the map  $\{\cdot,\cdot  \}  \colon L\times L\to  L $ is defined by $ \{x,y  \}=
    [\alpha(x),\alpha (y)],$
	for all $x,\, y\in L,$  then $(L,\{\cdot,\cdot  \},\alpha) $ is a symmetric Hom-Leibniz algebra, called $\alpha$-twist (or Yau twist) of $L$, and denoted by $L_{\alpha} $.
\end{prop}

Next, we consider the notions of right, left symmetric and symmetric (two-sided) Hom-Leibniz superalgebras.
\begin{defn}
	A triple  $(L,[\cdot,\cdot],\alpha) $ consisting of a superspace $ L $, an even  bilinear map $ [\cdot,\cdot]\, :L\times  L\rightarrow  L  $  and an even superspace homomorphism $ \alpha : L\rightarrow  L  $ \textup{(}linear map of parity degree $0\in\mathbb{Z}_2$\textup{)} is called	
	\begin{enumerate}[label=\upshape{(\roman*)},leftmargin=30pt]
		\item left Hom-Leibniz superalgebra if it satisfies for all $x, y, z \in L_{0}\cup L_1$,
		\begin{equation*}
		\left[\alpha(x),[y,z]\right]=	\left[ [x,y],\alpha(z)\right] + (-1)^{|x||y|} \left[\alpha(y), [x,z]\right],
		\end{equation*}
		\item right Hom-Leibniz superalgebra  if it satisfies
		%(see \cite{sami}) if satisfies
		\begin{equation*}
		\left[\alpha(x),[y,z]\right]=	\left[ [x,y],\alpha(z)\right] -(-1)^{|y||z|} \left[ [x,z],\alpha(y)\right],
		\end{equation*}
\item  symmetric Hom-Leibniz superalgebra if it is a left and a right Hom-Leibniz superalgebra.	
\end{enumerate}	
\end{defn}

\begin{prop}
	A triple $(L,[\cdot,\cdot],\alpha) $ is a symmetric
	Hom-Leibniz superalgebra  if and only if for all $ \ x,\ y,\ z\in L_{0}\cup L_1, $
	\begin{align*}
	\left[\alpha(x),[y,z]\right]&=	\left[ [x,y],\alpha(z)\right] + (-1)^{|x||y|} \left[\alpha(y), [x,z]\right],\\
	\left[\alpha(y), [x,z]\right]&=- (-1)^{\left( |x|+|z|\right) |y|} \left[ [x,z],\alpha(y)\right].
	\end{align*}
\end{prop}
\begin{example}
	Let $L=L_{0}\oplus L_{1}$  be a $3$-dimensional superspace, where $L_{0}$ is generated by
$e_1$, $e_2$ and $L_{1}$ is generated by $e_3$. The product is given by
	\begin{gather*}
	[e_1,e_1] =ae_1+xe_2, \quad
	[e_1,e_2]=[e_2,e_1]=-\frac{a}{x}[e_1,e_1],\quad
	[e_2,e_2]=\left(\frac{a}{x}\right)^{2}[e_1,e_1],\\
	[e_3,e_3] =\frac{d}{x}  [e_1,e_1],
\quad [e_1,e_3] =[e_3,e_1]=[e_3,e_2] =[e_2,e_3] =0  .
	\end{gather*}
Consider the homomorphism  $\alpha \colon L\to  L $ with the matrix $ \begin{pmatrix}
	-1&0&0\\
	0&1&0\\
	0&0&\mu
	\end{pmatrix} $ in the basis $ (e_1,e_2,e_3)$. Then $(L,[\cdot,\cdot],\alpha) $ is a symmetric
	Hom-Leibniz superalgebra.
\end{example}

%%%%%%%%%%%%%%%%%%%%%%%%%%%%%%%%%%%%%%%%
\section{Centroids and derivations of Hom-Leibniz superalgebras}
\label{sec:centroidsderivhomLeibnizsuperalgs}
In this section we consider centroids and twisted derivations of Hom-Leibniz superalgebras.
The concept of centroids and derivation of Leibniz algebras  is introduced in  \cite{SaidHusainRakhimovBasri2017:centrdsderlowdimleibnizalgs}. Left Leibniz superalgebras, originally were introduced by Dzhumadil’daev in \cite{Dzhumadildaev2000:CohomcolLeibnalgspresimplic}, can be seen as a direct generalization of Leibniz algebras. The left Hom-Leibniz superalgebras  were recently considered in \cite{PouyeJRM2021:HomLeibnizsuperinvnondbilfrms}.

Let $L$ be a superspace and $\alpha$ an even linear map on $L$.
Then $$\Omega = \{ u\in End(L)\mid u\circ\alpha=\alpha \circ u \}$$ is a linear subspace of $ End(L) $.
If the map $\tilde{\alpha} \colon  \Omega  \to     \Omega$ is defined by
$\tilde{\alpha}(u) =\alpha \circ u, $ then $ (\Omega,[\cdot,\cdot ]' ) $ is a Lie superalgebra, and $ (\Omega,[\cdot,\cdot ]', \tilde{\alpha})  $ ) is a Hom-Lie superalgebra with the bilinear bracket
defined by $ [u,v]'=uv-(-1)^{|u||v|} vu $ for all $ u,\, v \in \Omega. $	
The Hom-Lie algebras $ (\Omega,[\cdot,\cdot ]', \tilde{\alpha})  $ are not necessarily multiplicative since the linear map $\alpha$ does not need to be a homomorphism of the algebra structure defined by the bracket $ [\cdot,\cdot]'$.

\begin{defn}
For any non-negative integer $k\geq 0$, an $\alpha^k$-derivation of a  Hom-Leibniz superalgebra $(L,[\cdot,\cdot],\alpha) $ is a  homogeneous linear map  $D\in \Omega$ satisfying
\[ D([x,y])=[D(x),\alpha^{k}(y)]+(-1)^{|D||x|}[\alpha^{k}(x),D(y)] \]
for all $ \ x,\ y,\ z\in L_{0}\cup L_1. $
The set of all $\alpha^k$-derivations of a Hom-Leibniz superalgebra $L$ for all non-negative integers $k\geq 0$ is denoted by $\displaystyle Der(L)=\bigoplus_{k\geq0}Der_{\alpha^{k}}(L)  $.
\end{defn}
We will refer sometimes to the elements of $Der(L)$ as derivations of Hom-Leibniz superalgebra $L$ slightly abusing terminology for the convenience of the exposition.
\begin{prop}
Let  $(L,[\cdot,\cdot ],\alpha)$ be a left  (resp. right) Hom-Leibniz   superalgebra. For any $ a\in L $ satisfying $ \alpha(a)=a, $ define $ ad_{k}(a)\in End(L) $ and $ Ad_{k}(a)\in End(L) $ for all
$x\in L$ by
\begin{align*}
ad_{k}(a)(x) =[a,\alpha^{k}(x)],  \quad
Ad_{k}(a)(x) =(-1)^{|a||x|}[\alpha^{k}(x),a].
\end{align*}	
Then $ ad_{k}(a) $ \textup{(}resp. $ Ad_{k}(a) $\textup{)} is an $ \alpha^{k+1} $-derivation of the
left \textup{(}resp. right\textup{)} Hom-Leibniz algebra $ L $.
\end{prop}

\begin{defn}
Let  $(L,[\cdot,\cdot ],\alpha)$ be a Hom-Leibniz   superalgebra. Then the 	$\alpha^k$-centroid of $L$ denoted as $\Gamma_{\alpha^{k}}(L)$ is defined by \[\Gamma_{\alpha^{k}}(L)=\left\lbrace d\in  \Omega\mid  d([x,y])=[d(x),\alpha^{k}(y)] = (-1)^{|d||x|}[\alpha^{k}(x),d(y)],\,\forall  \ x,\ y\in L_{0}\cup L_1  \right\rbrace . \]
Denote by $\displaystyle \Gamma(L)=\bigoplus_{k\geq0}\Gamma_{\alpha^{k}}(L) $ the centroid of $\mathcal{L}$.
\end{defn}

\begin{defn}
	Let  $(L,[\cdot,\cdot ],\alpha)$ be a Hom-Leibniz   superalgebra. Then the 	$\alpha^k$-centroid of $L$ denoted by $C_{\alpha^{k}}(L)$ is defined by \[C_{\alpha^{k}}(L)=\left\lbrace d\in  \Omega\mid  d([x,y])=[d(x),\alpha^{k}(y)] = (-1)^{|d||x|}[\alpha^{k}(x),d(y)],\,\forall  \ x,y\in L_{0}\cup L_1                  \right\rbrace . \]
Denote by $\displaystyle C(L)=\bigoplus_{k\geq0}C_{\alpha^{k}}(L) $ the centroid of $L$.
\end{defn}

\begin{defn}
	Let  $(L,[\cdot,\cdot ],\alpha)$ be a Hom-Leibniz   superalgebra and $d\in End(L)$.
	Then $d$ is called $\alpha^k$-central derivation, if $d\in \Omega$ and \[ d([x,y])=[d(x),\alpha^{k}(y)]=(-1)^{|d||x|}[\alpha^{k}(x),d(y)]=0.\]
The set of all  $\alpha^k$-central derivations of a Hom-Leibniz superalgebra $L$ for all non-negative integers $k\geq 0$ is denoted by $ \displaystyle ZDer(L)=\bigoplus_{k\geq0}ZDer_{\alpha^{k}}(L) $.	
\end{defn}

In the rest of the section, we study the structure of the centroids and derivations of  Hom-Leibniz superalgebras.
\begin{lem}
	If $L$ is a  Hom-Leibniz superalgebra, $d\in Der_{\alpha^{k}}(L)  $ and
	$\Phi\in C_{\alpha^{l}}(L)  $, then
	\begin{enumerate}[label=\upshape{(\roman*)},leftmargin=30pt]
		\item $\Phi\circ d$ is an $\alpha^{k+l}$-derivation of $L$.
		\item $[\Phi,d]$ is $\alpha^{k+l}$-centroid of $L$.
		\item $ d\circ \Phi$ is an $\alpha^{k+l}$-centroid if and only if $\Phi\circ d$ is a $\alpha^{k+l}$-central derivation.
		\item $ d\circ \Phi$ is an $\alpha^{k+l}$-derivation only if only $[d,\Phi ]$ is a $\alpha^{k+l}$-central derivation.
	\end{enumerate}
\end{lem}

The proof of the following theorem is similar to the case of Leibniz algebra given in \cite{SaidHusainRakhimovBasri2017:centrdsderlowdimleibnizalgs}.
\begin{thm}
	If $L$ is a  Hom-Leibniz superalgebra, then \[ZDer_{\alpha^{k}}(L) =Der_{\alpha^{k}}(L)\cap C_{\alpha^{k}}(L).\]
\end{thm}

%\begin{prop}
%	Let  $(L,[\cdot,\cdot ],\alpha,\beta)$ be a left  BiHom-Leibniz   algebra and $ a\in L $.
%	Define respectively $ D,\quad D',\quad D'' $
%	by
%	\begin{equation*}
%	D(x)=[\alpha\beta(a),\alpha^{k}\beta^{l}(x)],\quad  D'(x)=[\beta(a),\alpha^{k}\beta^{l}(x)],\quad   D''(x)=[\alpha(a),\alpha^{k}\beta^{l}(x)].
%	\end{equation*}	
%	Then
%\begin{enumerate}[(i)]
%	\item \quad  $ \alpha\circ D'=D\circ\alpha, $ $ \beta\circ D''=D\circ\beta $
%	\item \quad$ 	[D(x),\alpha^{k}\beta^{l}(y)]+[\alpha^{k}\beta^{l}(x),D'(y)]=D''([x,y]), $
%	\item \quad If $ \alpha^{2}(a) =\alpha(a)=\beta(a)$, $ D $ is a $\alpha^{k}\beta^{l+1}  $-derivation   of the
%	left  BiHom-Leibniz algebra $ L $.
%	\end{enumerate}
%\end{prop}
%\date{ 18/07/2018}\\

\begin{prop}
	If $(L,[\cdot,\cdot ],\alpha)$ is a Hom-Leibniz superalgebra, then
	\begin{enumerate}[label=\upshape{(\roman*)},leftmargin=30pt]
		\item $ ad(L)\subseteq  Der(L) \subseteq \Omega$, where $ad(L)  $ denotes the superalgebra of inner derivations of $ L $.
		\item  $ ad(L) $, 	$ Der(L) $  and $ C(L) $ are Lie \textup{(}resp. Hom-Lie\textup{)} subsuperalgebras of $ (\Omega,[\cdot,\cdot ]' ) $ \textup{(}resp. $(\Omega,[\cdot,\cdot ]', \tilde{\alpha})$\textup{)}.
	\end{enumerate}
	
\end{prop}
%%%%%%%%%%%%%%%%%%%%%%%%%%%%%%%%%%%%%%%%%%%%%%%%%%
%%%%%%%%%%%%%%%%%%%%%%%%%%%%%%%%
%%%%%%%%%%%%%%%%%%%%%%%%%%%%%%%%%%%%%%%%
%\begin{prop}	
%We assume that $ \alpha^2=\alpha $ 	algebra and define the
%following linear space $ \Omega$ of $ End(L) $ consisting of linear maps $u$ on $ L $ as follows
%\[ \Omega =\{ u\in End(L)/ u\alpha=\alpha u    \}.\]
%Then $ \Omega $	 is a
%symmetric
%Hom-Lie
%Leibniz
% superalgebra with the bracket  \[[u,v]=uv-(-1)^{|u||v|} vu\]
%and the map \[ \sigma\colon  \Omega  \to     \Omega;\quad \sigma(u) =\alpha u.       \]
%\end{prop}	
%\begin{prop}\cite{superleibniz}
%	Let $V $ be a Hom-Lie superalgebra, then  the bracket
%	\[ [x\otimes y,a\otimes b]=[[x,y],\alpha(a)]\otimes b+(-1)^{\arrowvert a\arrowvert(\arrowvert x\arrowvert+\arrowvert y\arrowvert)}a\otimes[[x,y],\alpha(b)] \] defines  a right
%	Hom-Leibniz superalgebra structure on the linear space $ V\otimes V $.	
%\end{prop}
%
%\subsection{Hom-Leibniz colour algebras}

%\subsection{        }

%\subsection{ Hom Poisson Algebras      }
%\subsection{ }
%\subsection{ }
%\textbf{Voir les autres types et chercher un relation Leibniz et symmetric Leibniz}
%%%%%%%%%%%%%%%%%%%%%%%%%B2B2B2B2B2B2222222222222222222222222222222222222222222

%%%%%%%%%%%%%%%%%%%%%%%%%%%%%%%%ù
\section{Classification of multiplicative 2-dimensional Hom-Leibniz algebras}
\label{sec:clasmult2dimhomLeibnizalgs}

First of all, note that the classifications of two and three-dimensional complex Leibniz algebras were studied in \cite{LodayEM1993:versnoncomalgLie,BenkartNeherJA2006:CantroidExtAffinerootgrad}.
In this section, the classification of the $2$-dimensional left Hom-Leibniz algebras is obtained, and for each isomorphism class it is indicated whether the Hom-Leibniz algebras from this class are symmetric Hom-Leibniz algebras or not.

\begin{prop}
	Every  $ 2 $-dimensional left  Hom-Leibniz algebra is isomorphic to one of the following nonisomorphic Hom-Leibniz algebras, with
    each  algebra denoted by $ L_{j}^{i} $ where $ i $ is related to the linear map $ \alpha $ and $ j $ being the number in the list. \\
	%%%%%%%%%%%%%%%%%%%%%

	%%%%%%%%%%%%%%%%%%%%%%%%%%%%%%%%%%%%%%%%%%%%%%%%%%%%%%%%%%%
	%alg1 y1\neq 0 et x2y1+z2^2\neq 0
	%%%%%%%%%%%%%%%%%%%%%%%%%%%%%%%%%%%%%%%%%%%%%%%%%%%%
	\hrule
	\begin{align*}
	L_{1}^{1}:[e_{1},e_{1}]&=0,&     [e_{1},e_{2}]&=0,
	&
	[e_2,e_{1}]&=e_{1}
	,&       [e_2,e_{2}]&=e_1,&\\
	\alpha(e_{1})&=e_{1},& \alpha(e_{2})&=e_{2}.&&&&   	
	\end{align*} \hspace{1,5cm}	
$ L_{1}^{1} $ is not symmetric \\
	%%%%%%%%%%%%%%%%%%%%%
	\hrule
	%%%%%%%%%%%%%%%%%%%%%%%%%%%%%%%%%%%%%%%%%%%%%%%%%%%%%%%%%%%
	%alg1 y1\neq 0 et x2y1+z2^2= 0
	%%%%%%%%%%%%%%%%%%%%%%%%%%%%%%%%%%%%%%%%%%%%%%%%%%%%	
	\begin{align*}
	L_{2}^{1}:[e_{1},e_{1}]&=0,&     [e_{1},e_{2}]&=0,
	&
	[e_2,e_{1}]&=0
	,&       [e_2,e_{2}]&=e_1,&\\
	\alpha(e_{1})&=e_{1},& \alpha(e_{2})&=e_{2}.&&&&   	
	\end{align*} \hspace{1,5cm}
	$ L_{2}^{1} $ is  symmetric \\
	%%%%%%%%%%%%%%%%%%%%%Alg2 $ \alpha=id $	
	\hrule
	%%%%%%%%%%%%%%%%%%%%%%%%%%%%%%%%%%%%%%%%%%%%%%%%%%%%%%%%%%%
	%alg1 y1= 0
	%%%%%%%%%%%%%%%%%%%%%%%%%%%%%%%%%%%%%%%%%%%%%%%%%%%%	
	\begin{align*}
	L_{3}^{1}:[e_{1},e_{1}]&=0,&     [e_{1},e_{2}]&=0,
	&
	[e_2,e_{1}]&=e_{1}
	,&       [e_2,e_{2}]&=0,&\\
	\alpha(e_{1})&=e_{1},& \alpha(e_{2})&=e_{2}.&&&&   	
	\end{align*} \hspace{1,5cm}
	$ L_{3}^{1} $ is  symmetric	\\
	%%%%%%%%%%%%%%%%%%%%%Alg6 $ \alpha=id $, $ y1\neq 0 $	
	\hrule	
	\begin{align*}
	L_{4}^{1}:[e_{1},e_{1}]&=0,&     [e_{1},e_{2}]&=e_1,
	&
	[e_2,e_{1}]&=-e_{1}
	,&       [e_2,e_{2}]&=0,&\\
	\alpha(e_{1})&=e_{1},& \alpha(e_{2})&=e_{2}.&&&&   	
	\end{align*} \hspace{1,5cm}
	$L_{4}^{1}  $ is symmetric	\\
	%%%%%%%%%%%%%%%%%%%%%Alg1. $ a=0,b=1 $,
	\hrule	
	\begin{align*}
	L_{1}^{2}:[e_{1},e_{1}]&=0,&     [e_{1},e_{2}]&=0,
	&
	[e_2,e_{1}]&=e_{1}
	,&       [e_2,e_{2}]&=0,&\\
	\alpha(e_{1})&=0,& \alpha(e_{2})&=e_{2}.&&&&   	
	\end{align*} \hspace{1,5cm}
	$L_{1}^{2}  $ is not symmetric	\\
	\hrule	
	\begin{align*}
	L_{2}^{2}:[e_{1},e_{1}]&=e_1,&     [e_{1},e_{2}]&=0,
	&
	[e_2,e_{1}]&=0
	,&       [e_2,e_{2}]&=0,&\\
	\alpha(e_{1})&=0,& \alpha(e_{2})&=e_{2}.&&&&   	
	\end{align*} \hspace{1,5cm}
	$L_{2}^{2}  $ is  symmetric \\
	\hrule	
	\begin{align*}
	L_{3}^{2}:[e_{1},e_{1}]&=0,&     [e_{1},e_{2}]&=e_1,
	&
	[e_2,e_{1}]&=-e_{1}
	,&       [e_2,e_{2}]&=0,&\\
	\alpha(e_{1})&=0,& \alpha(e_{2})&=e_{1}.&&&&   	
	\end{align*} \hspace{1,5cm}
	$L_{3}^{2}  $ is  symmetric	\\
	\hrule	
	\begin{align*}
	L_{1}^{3}:[e_{1},e_{1}]&=e_1,&     [e_{1},e_{2}]&=0,
	&
	[e_2,e_{1}]&=0
	,&       [e_2,e_{2}]&=0,&\\
	\alpha(e_{1})&=0,& \alpha(e_{2})&=be_{2}.&&&&   	
	\end{align*} \hspace{1,5cm}
	$L_{1}^{3}  $ is  symmetric \\
	\hrule	
	\begin{align*}
	L_{2}^{3}:[e_{1},e_{1}]&=0,&     [e_{1},e_{2}]&=0,
	&
	[e_2,e_{1}]&=e_1
	,&       [e_2,e_{2}]&=0,&\\
	\alpha(e_{1})&=0,& \alpha(e_{2})&=be_{2}.&&&&   	
	\end{align*} \hspace{1,5cm}
	$L_{2}^{3}  $ is  not symmetric \\
	\hrule	
	\begin{align*}
	L_{3}^{3}:[e_{1},e_{1}]&=0,&     [e_{1},e_{2}]&=e_1,
	&
	[e_2,e_{1}]&=-e_1
	,&       [e_2,e_{2}]&=0,&\\
	\alpha(e_{1})&=0,& \alpha(e_{2})&=be_{2}.&&&&   	
	\end{align*} \hspace{1,5cm}
	$L_{3}^{3}  $ is   symmetric \\
	\hrule	
	\begin{align*}
	L_{1}^{4}:[e_{1},e_{1}]&=0,&     [e_{1},e_{2}]&=0,
	&
	[e_2,e_{1}]&=e_1
	,&       [e_2,e_{2}]&=0,&\\
	\alpha(e_{1})&=ae_1 (a\notin\{0,1\}),& \alpha(e_{2})&=e_{2}.&&&&   	
	\end{align*} \hspace{1,5cm}
	$L_{1}^{4}  $ is not  symmetric \\
	\hrule	
	\begin{align*}
	L_{2}^{4}:[e_{1},e_{1}]&=0,&     [e_{1},e_{2}]&=-e_1,
	&
	[e_2,e_{1}]&=e_1
	,&       [e_2,e_{2}]&=0,&\\
	\alpha(e_{1})&=ae_1 (a\notin\{0,1\}),& \alpha(e_{2})&=e_{2}.&&&&   	
	\end{align*} \hspace{1,5cm}
	$L_{2}^{4}  $ is not  symmetric \\
	\hrule	
	\begin{align*}
	L_{1}^{5}:[e_{1},e_{1}]&=0,&     [e_{1},e_{2}]&=0,
	&
	[e_2,e_{1}]&=0
	,&       [e_2,e_{2}]&=e_1,&\\
	\alpha(e_{1})&=a^{2}e_1 (a\notin\{0,1\}),& \alpha(e_{2})&=ae_{2}.&&&&   	
	\end{align*} \hspace{1,5cm}
	$L_{1}^{5}  $ is   symmetric \\
\hrule	
\begin{align*}
	L_{1}^{6}:[e_{1},e_{1}]&=0,&     [e_{1},e_{2}]&=e_{1},
	&
	[e_2,e_{1}]&=y_1e_1
	,&       [e_2,e_{2}]&=z_1e_1,&\\
	\alpha(e_{1})&=0 ,& \alpha(e_{2})&=e_1.&&&&   	
\end{align*} \hspace{1,5cm}
$L_{1}^{6}  $ is   symmetric	\\
\hrule	
\begin{align*}
L_{2}^{6}:[e_{1},e_{1}]&=0,&     [e_{1},e_{2}]&=0,
&
[e_2,e_{1}]&=e_1
,&       [e_2,e_{2}]&=z_1e_1,&\\
\alpha(e_{1})&=0 ,& \alpha(e_{2})&=e_1.&&&&   	
\end{align*} \hspace{1,5cm}
$L_{2}^{6}  $ is   symmetric	\\
\hrule	
\begin{align*}
L_{3}^{6}:[e_{1},e_{1}]&=0,&     [e_{1},e_{2}]&=0,
&
[e_2,e_{1}]&=0
,&       [e_2,e_{2}]&=e_1,&\\
\alpha(e_{1})&=0 ,& \alpha(e_{2})&=e_1.&&&&   	
\end{align*} \hspace{1,5cm}
$L_{3}^{6}  $ is   symmetric	\\
\hrule	
\begin{align*}
L_{1}^{7}:[e_{1},e_{1}]&=0,&     [e_{1},e_{2}]&=0,
&
[e_2,e_{1}]&=0
,&       [e_2,e_{2}]&=e_1,&\\
\alpha(e_{1})&=e_{1} ,& \alpha(e_{2})&=e_1+e_2.&&&&   	
\end{align*} \hspace{1,5cm}	
$L_{1}^{7}  $ is   symmetric \\	
\hrule
\end{prop}	

\section{Centroids and derivations of  2-dimensional multiplicative Hom-Leibniz algebras}
\label{sec:centroidsderiv2dimmulthomLeibnizalgs}
In this section we focus on the study of centroids and derivation of multiplicative Hom-Leibniz algebras of dimension $2$ over $\mathbb{C}$. Note that if the odd component of  Hom-Leibniz superalgebra is zero, it can be considered as Hom-Leibniz algebra.

Next, we introduce a parametric extension of the notion of $ \alpha^k $-derivations of Hom-Leibniz algebras and study it in more details with help of the computer computations.
	
\begin{defn} Let $ (L,[\cdot,\cdot],\alpha) $ be a Hom–Leibniz algebra and $\lambda,\mu,\gamma $
		elements of $\C$.
		A linear map $d\in \Omega$ is a generalized 	$ \alpha^k $-derivation or a	
		$ (\lambda,\mu,\gamma)$-$ \alpha^k$-derivation of $L$ if for all $x,y\in L$ we have\[ \lambda\,	d([x,y])=\mu\, [d(x),\alpha^{k}(y)]+ \gamma\,[\alpha^{k}(x),d(y)]. \]
		%	\textbf{ATTENTION THE SUM IS NOT  necessarily direct}
		We denote the set of all $ (\lambda,\mu,\gamma)$-$ \alpha^k $-derivations by $\displaystyle  Der_{\alpha^k}^{(\lambda,\mu,\gamma)}(L)$.
		%and $\displaystyle Der^{(\lambda,\mu,\gamma)}(L) $ the vector space spanned by  $ \{d\in Der_{\alpha^k\beta^l}^{(\lambda,\mu,\gamma)}(L)\mid k,l\in \N\} $.		
		%the set of generalized deivation of $ L $.
\end{defn}

Clearly $\displaystyle  Der_{\alpha^k}^{(1,1,1)}(L)=Der_{\alpha^k}(L)$ and $\displaystyle  Der_{\alpha^k}^{(1,1,0)}(L)=\Gamma_{\alpha^k}(L)$. Let  $(L,[\cdot,\cdot ],\alpha)$ be a
$n$-dimensional multiplicative left Hom-Leibniz algebra. Let $\displaystyle \alpha^{r}(e_{j})=\sum_{k=1}^{n}m_{kj}e_{k} $.
An element $ d $ of $ Der^{(\delta,\mu,\gamma)}_{\alpha^{r}}(L)$, being a linear transformation of the linear space $ L$, is represented in a matrix form $ (d_{ij})_{1\leq i,j\leq n} $, that is $\displaystyle d(e_{j})=\sum_{k=1}^{n}d_{kj}e_{k} $, for $ j=1,\cdots,n $. According to the definition of the $ (\delta,\mu,\gamma) $-$ \alpha^{r} $-derivation the entries $ d_{ij} $ of the matrix $ (d_{ij})_{1\leq i,j\leq n} $
must satisfy the following systems of equations:
\begin{align*}
\displaystyle \sum_{k=1}^{n} d_{ik}a_{kj} &= \sum_{k=1}^{n} a_{ik}d_{kj},\\
\displaystyle \delta\sum_{k=1}^{n} c_{ij}^{k}d_{sk}&-\mu \sum_{k=1}^{n} \sum_{l=1}^{n}d_{ki}m_{lj}c_{kl}^{s} -\gamma   \sum_{k=1}^{n} \sum_{l=1}^{n}d_{lj} m_{ki}c_{kl}^{s}=0,
\end{align*}
where $ (a_{ij})_{1\leq i,j\leq n} $ is the matrix of $ \alpha $ and
$( c_{ij}^{k}) $ are the structure constants of $ L $.
%\subsection{Centroids  of  $ 2 $-dimensional multiplicative BiHom-Lie algebras}	
First, let us give the following definitions.
\begin{defn}
	A left Hom-Leibniz algebra  is called characteristically nilpotent if the Lie algebra $ Der_{\alpha^{0}}(L) $ is nilpotent.
\end{defn}
Henceforth, the property of being characteristically nilpotent is abbreviated by CN.

\begin{defn}
	Let $ L $ be an indecomposable left Hom-Leibniz algebra. We say that $ L $
	is small if $\varGamma_{\alpha^{0}}(L)  $ is generated by central derivation and the scalars. 	
	The centroid of a decomposable 	BiHom-Lie algebra is small if the centroids of each indecomposable factor is small.	
\end{defn}
Now we apply the algorithms mentioned in the previous paragraph to find centroid and derivations of
$ 2 $-dimensional complex left Hom-Leibniz
algebras. In this study we make use of the classification results from the previous section.
The results are given in the following theorem. Moreover, we provide the type of  $ \varGamma_{\alpha^0}(L_{i}^{j}) $ and $ Der_{\alpha^r}(L_{i}^{j})  $. %if $ (r,l)=(0,0) $.
\begin{align*}
L_{1}^{1}:[e_{1},e_{1}]&=0,&     [e_{1},e_{2}]&=0,
&
[e_2,e_{1}]&= e_{1}
,&       [e_2,e_{2}]&=e_{1},&\\
\alpha(e_{1})&=e_{1},& \alpha(e_{2})&=e_{2}&&& &&   	
\end{align*}
\begin{center}
	\begin{tabular}{|c|c|c|c|c|}	
		\hline
		$ \alpha^r $ &$ \varGamma_{\alpha^r}(L_{1}^{1}) $& Type of $ \varGamma_{\alpha^0}(L_{1}^{1}) $    &$ Der_{\alpha^r}(L_{1}^{1}) $&CN\\
		\hline
		$ r\in \N$&$ \begin{pmatrix}
		c_1&c_2\\
		0&c_1
		\end{pmatrix} $	&not  small   &$ \begin{pmatrix}
		d_1&d_1\\
		0&0
		\end{pmatrix} $& Yes\\		
		\hline
	\end{tabular}
\end{center}
%%%%%%%%%%%%%%%%%%%%%	
\begin{align*}
L_{2}^{1}:[e_{1},e_{1}]&=0,&     [e_{1},e_{2}]&=0,
&
[e_2,e_{1}]&= 0
,&       [e_2,e_{2}]&=e_{1},&\\
\alpha(e_{1})&=e_{1},& \alpha(e_{2})&=e_{2}&&& &&   	
\end{align*}
\begin{center}
	\begin{tabular}{|c|c|c|c|c|}	
		\hline
		$ \alpha^r $ &$ \varGamma_{\alpha^r}(L_{1}^{1}) $& Type of $ \varGamma_{\alpha^0}(L_{2}^{1}) $    &$ Der_{\alpha^r}(L_{1}^{1}) $&CN\\
		\hline
		$ r\in \N$&$ \begin{pmatrix}
		c_1&c_2\\
		0&c_1
		\end{pmatrix} $	&  small   &$ \begin{pmatrix}
		d_1&d_2\\
		0&\frac{d_1}{2}
		\end{pmatrix} $& Yes\\[0.35cm]		
		\hline
	\end{tabular}
\end{center}
%%%%%%%%%%%%%%%%%%%%%
\begin{align*}
L_{3}^{1}:[e_{1},e_{1}]&=0,&     [e_{1},e_{2}]&=0,
&
[e_2,e_{1}]&= e_1
,&       [e_2,e_{2}]&=0,&\\
\alpha(e_{1})&=e_{1},& \alpha(e_{2})&=e_{2}&&& &&   	
\end{align*} 	
\begin{center}
	\begin{tabular}{|c|c|c|c|c|}	
		\hline
		$ \alpha^r $ &$ \varGamma_{\alpha^r}(L_{3}^{1}) $& Type of $ \varGamma_{\alpha^0}(L_{2}^{1}) $    &$ Der_{\alpha^r}(L_{3}^{1}) $&CN\\
		\hline
		$ r\in \N$&$ \begin{pmatrix}
		c_1&c_2\\
		0&c_1
		\end{pmatrix} $	& not small   &$ \begin{pmatrix}
		d_1&0\\
		0&0
		\end{pmatrix} $& Yes\\		
		\hline
	\end{tabular}
\end{center}
%%%%%%%%%%%%%%%%%%%%%	
%%%%%%%%%%%%%%%%%%%%%
\begin{align*}
L_{4}^{1}:[e_{1},e_{1}]&=0,&     [e_{1},e_{2}]&=e_1,
&
[e_2,e_{1}]&= -e_1
,&       [e_2,e_{2}]&=0,&\\
\alpha(e_{1})&=e_{1},& \alpha(e_{2})&=e_{2}&&& &&   	
\end{align*} 	
\begin{center}
	\begin{tabular}{|c|c|c|c|c|}	
		\hline
		$ \alpha^r $ &$ \varGamma_{\alpha^r}(L_{4}^{1}) $& Type of $ \varGamma_{\alpha^0}(L_{4}^{1}) $    &$ Der_{\alpha^r}(L_{4}^{1}) $&CN\\
		\hline
		$ r\in \N$&$ \begin{pmatrix}
		c_1&0\\
		0&c_1
		\end{pmatrix} $	&  small   &$ \begin{pmatrix}
		d_1&d_2\\
		0&0
		\end{pmatrix} $& No\\		
		\hline
	\end{tabular}
\end{center}
%%%%%%%%%%%%%%%%%%%%%	
%%%%%%%%%%%%%%%%%%%%%
\begin{align*}
L_{1}^{2}:[e_{1},e_{1}]&=0,&     [e_{1},e_{2}]&=0,
&
[e_2,e_{1}]&= e_1
,&       [e_2,e_{2}]&=0,&\\
\alpha(e_{1})&=0,& \alpha(e_{2})&=e_{2}&&& &&   	
\end{align*} 	
\begin{center}
	\begin{tabular}{|c|c|c|c|c|}	
		\hline
		$ \alpha^r $ &$ \varGamma_{\alpha^r}(L_{1}^{2}) $& Type of $ \varGamma_{\alpha^0}(L_{1}^{2}) $    &$ Der_{\alpha^r}(L_{1}^{2}) $&CN\\
		\hline
		$ r\in \N$&$ \begin{pmatrix}
		0&0\\
		0&c_2
		\end{pmatrix} $	& not small   &$ \begin{pmatrix}
		d_1&0\\
		0&d_2
		\end{pmatrix} $& Yes\\		
		\hline
	\end{tabular}
\end{center}
%%%%%%%%%%%%%%%%%%%%%
\begin{align*}
L_{2}^{2}:[e_{1},e_{1}]&=e_1,&     [e_{1},e_{2}]&=0,
&
[e_2,e_{1}]&= 0
,&       [e_2,e_{2}]&=0,&\\
\alpha(e_{1})&=0,& \alpha(e_{2})&=e_{2}&&& &&   	
\end{align*} 	
\begin{center}
	\begin{tabular}{|c|c|c|c|c|}	
		\hline
		$ \alpha^r $ &$ \varGamma_{\alpha^r}(L_{2}^{2}) $& Type of $ \varGamma_{\alpha^0}(L_{2}^{2}) $    &$ Der_{\alpha^r}(L_{2}^{2}) $&CN\\
		\hline
		$ r\in \N$&$ \begin{pmatrix}
		0&0\\
		0&c_2
		\end{pmatrix} $	&  small   &$ \begin{pmatrix}
		0&0\\
		0&d_2
		\end{pmatrix} $& Yes\\		
		\hline
	\end{tabular}
\end{center}
%%%%%%%%%%%%%%%%%%%%%
\begin{align*}
L_{3}^{2}:[e_{1},e_{1}]&=e_1,&     [e_{1},e_{2}]&=0,
&
[e_2,e_{1}]&= 0
,&       [e_2,e_{2}]&=0,&\\
\alpha(e_{1})&=0,& \alpha(e_{2})&=e_{2}&&& &&   	
\end{align*} 	
\begin{center}
	\begin{tabular}{|c|c|c|c|c|}	
		\hline
		$ \alpha^r $ &$ \varGamma_{\alpha^r}(L_{3}^{2}) $& Type of $ \varGamma_{\alpha^0}(L_{3}^{2}) $    &$ Der_{\alpha^r}(L_{3}^{2}) $&CN\\
		\hline
		$ r\in \N$&$ \begin{pmatrix}
		0&0\\
		0&c_2
		\end{pmatrix} $	& not small   &$ \begin{pmatrix}
		d_1&0\\
		0&d_2
		\end{pmatrix} $& Yes\\		
		\hline
	\end{tabular}
\end{center}
%%%%%%%%%%%%%%%%%%%%%
\begin{align*}
L_{1}^{3}:[e_{1},e_{1}]&=e_1,&     [e_{1},e_{2}]&=0,
&
[e_2,e_{1}]&= 0
,&       [e_2,e_{2}]&=0,&\\
\alpha(e_{1})&=0,& \alpha(e_{2})&=be_{2}&&& &&   	
\end{align*} 	
\begin{center}
	\begin{tabular}{|c|c|c|c|c|}	
		\hline
		$ \alpha^r $ &$ \varGamma_{\alpha^r}(L_{1}^{3}) $& Type of $ \varGamma_{\alpha^0}(L_{1}^{3}) $    &$ Der_{\alpha^r}(L_{1}^{3}) $&CN\\
		\hline
		$ r\in \N$&$ \begin{pmatrix}
		0&0\\
		0&c_2
		\end{pmatrix} $	&  small   &$ \begin{pmatrix}
		0&0\\
		0&d_2
		\end{pmatrix} $& Yes\\		
		\hline
	\end{tabular}
\end{center}
%%%%%%%%%%%%%%%%%%%%%
\begin{align*}
L_{2}^{3}:[e_{1},e_{1}]&=0,&     [e_{1},e_{2}]&=0,
&
[e_2,e_{1}]&=e_1
,&       [e_2,e_{2}]&=0,&\\
\alpha(e_{1})&=0,& \alpha(e_{2})&=be_{2}\, (b\neq 1)&&& &&   	
\end{align*} 	
\begin{center}
	\begin{tabular}{|c|c|c|c|c|}	
		\hline
		$ \alpha^r $ &$ \varGamma_{\alpha^r}(L_{2}^{3}) $& Type of $ \varGamma_{\alpha^0}(L_{2}^{3}) $    &$ Der_{\alpha^r}(L_{2}^{3}) $&CN\\
		\hline
		$ r\in \N$&$ \begin{pmatrix}
		0&0\\
		0&c_2
		\end{pmatrix} $	& not small   &$ \begin{pmatrix}
		0&0\\
		0&d_2
		\end{pmatrix} $& Yes\\		
		\hline
	\end{tabular}
\end{center}
%%%%%%%%%%%%%%%%%%%%%
\begin{align*}
L_{3}^{3}:[e_{1},e_{1}]&=0,&     [e_{1},e_{2}]&=e_1,
&
[e_2,e_{1}]&=-e_1
,&       [e_2,e_{2}]&=0,&\\
\alpha(e_{1})&=0,& \alpha(e_{2})&=be_{2}\, (b\neq 1)&&& &&   	
\end{align*} 	
\begin{center}
	\begin{tabular}{|c|c|c|c|c|}	
		\hline
		$ \alpha^r $ &$ \varGamma_{\alpha^r}(L_{3}^{3}) $& Type of $ \varGamma_{\alpha^0}(L_{3}^{3}) $    &$ Der_{\alpha^r}(L_{3}^{3}) $&CN\\
		\hline
		$ r\in \N$&$ \begin{pmatrix}
		0&0\\
		0&c_2
		\end{pmatrix} $	& not small   &$ \begin{pmatrix}
		0&0\\
		0&d_2
		\end{pmatrix} $& Yes\\		
		\hline
	\end{tabular}
\end{center}
%%%%%%%%%%%%%%%%%%%%%
\begin{align*}
L_{1}^{4}:[e_{1},e_{1}]&=0,&     [e_{1},e_{2}]&=0,
&
[e_2,e_{1}]&=e_1
,&       [e_2,e_{2}]&=0,&\\
\alpha(e_{1})&=ae_1,& \alpha(e_{2})&=e_{2}&&& &&   	
\end{align*} 	
\begin{center}
	\begin{tabular}{|c|c|c|c|c|}	
		\hline
		$ \alpha^r $ &$ \varGamma_{\alpha^r}(L_{1}^{4}) $& Type of $ \varGamma_{\alpha^0}(L_{1}^{4}) $    &$ Der_{\alpha^r}(L_{1}^{4}) $&CN\\
		\hline
		$ r\in \N$&
        \begin{minipage}{2cm} $
        \begin{pmatrix}
		c_1&0\\
		0&\frac{c_1}{a^r}
		\end{pmatrix} $	\end{minipage} &  small   &$ \begin{pmatrix}
		0&0\\
		0&d_2
		\end{pmatrix} $& Yes \\[0.35cm]		
		\hline
	\end{tabular}
\end{center}
%%%%%%%%%%%%%%%%%%%%%
\begin{align*}
L_{2}^{4}:[e_{1},e_{1}]&=0,&     [e_{1},e_{2}]&=-e_1,
&
[e_2,e_{1}]&=e_1
,&       [e_2,e_{2}]&=0,&\\
\alpha(e_{1})&=ae_1,& \alpha(e_{2})&=e_{2}&&& &&   	
\end{align*} 	
\begin{center}
	\begin{tabular}{|c|c|c|c|c|}	
		\hline
		$ \alpha^r $ &$ \varGamma_{\alpha^r}(L_{2}^{4}) $& Type of $ \varGamma_{\alpha^0}(L_{2}^{4}) $    &$ Der_{\alpha^r}(L_{2}^{4}) $&CN\\
		\hline
		$ r\in \N$&$ \begin{pmatrix}
		c_1&0\\
		0&\frac{c_1}{a^r}
		\end{pmatrix} $	&  small   &$ \begin{pmatrix}
	d_1	&0\\
		0&0
		\end{pmatrix} $& Yes\\[0.35cm]		
		\hline
	\end{tabular}
\end{center}
%%%%%%%%%%%%%%%%%%%%%
\begin{align*}
L_{1}^{5}:[e_{1},e_{1}]&=0,&     [e_{1},e_{2}]&=0,
&
[e_2,e_{1}]&=0
,&       [e_2,e_{2}]&=e_1,&\\
\alpha(e_{1})&=b^2e_1,& \alpha(e_{2})&=be_{2}&&& &&   	
\end{align*} 	
\begin{center}
	\begin{tabular}{|c|c|c|c|c|}	
		\hline
		$ \alpha^r $ &$ \varGamma_{\alpha^r}(L_{1}^{5}) $& Type of $ \varGamma_{\alpha^0}(L_{1}^{5}) $    &$ Der_{\alpha^r}(L_{1}^{5}) $&CN\\
		\hline
		$ r\in \N$&$ \begin{pmatrix}
		c_1&0\\
		0&\frac{c_1}{b^r}
		\end{pmatrix} $	 &  small   &$ \begin{pmatrix}
		0	&0\\
		0&0
		\end{pmatrix} $ & Yes \\[0.35cm] 		
		\hline
	\end{tabular}
\end{center}
%%%%%%%%%%%%%%%%%%%%%
\begin{align*}
L_{1}^{6}:[e_{1},e_{1}]&=0,&     [e_{1},e_{2}]&=e_1,
&
[e_2,e_{1}]&=z_1e_1
,&       [e_2,e_{2}]&=t_1e_1,&\\
\alpha(e_{1})&=0,& \alpha(e_{2})&=e_{2} &&& &&   	
\end{align*} 	
\begin{center}
	\begin{tabular}{|c|c|c|c|c|c|}	
		\hline
		$ \alpha^r $& &$ \varGamma_{\alpha^r}(L_{1}^{6}) $& Type of $ \varGamma_{\alpha^0}(L_{1}^{6}) $    &$ Der_{\alpha^r}(L_{1}^{6}) $&CN\\
		\hline
		$ r=0$&$ z_1=-1 $&$ \begin{pmatrix}
		c_1&0\\
		0&c_1
		\end{pmatrix} $	&  small   &$ \begin{pmatrix}
		0	&d_2\\
		0&0
		\end{pmatrix} $& Yes\\
	\hline
$ r=0$&$ z_1\neq-1 $&$ \begin{pmatrix}
c_1&0\\
0&c_1
\end{pmatrix} $	&  small   &$ \begin{pmatrix}
0	&0\\
0&0
\end{pmatrix} $& Yes\\				
		%\hline
%$ \alpha^r $& &$ \varGamma_{\alpha^r}(L_{1}^{6}) $& Type of $ \varGamma_{\alpha^0}(L_{1}^{6}) $    &$ Der_{\alpha^r}(L_{1}^{6}) $&CN\\
\hline
$ r=1$&&$ \begin{pmatrix}
0&c_2\\
0&0
\end{pmatrix} $	&   &$ \begin{pmatrix}
0	&d_2\\
0&0
\end{pmatrix} $& \\		
\hline
$ r>1$&&$ \begin{pmatrix}
0&c_2\\
0&0
\end{pmatrix} $	&  &$ \begin{pmatrix}
0	&d_2\\
0&0
\end{pmatrix} $& \\		
\hline				
	\end{tabular}
\end{center}
%%%%%%%%%%%%%%%%%%%%%
\begin{align*}
L_{2}^{6}:[e_{1},e_{1}]&=0,&     [e_{1},e_{2}]&=0,
&
[e_2,e_{1}]&=e_{1}
,&       [e_2,e_{2}]&=t_1e_1,&\\
\alpha(e_{1})&=0,& \alpha(e_{2})&=e_{1}&&& &&   	
\end{align*} 	
\begin{center}
	\begin{tabular}{|c|c|c|c|c|}	
		\hline
		$ \alpha^r $ &$ \varGamma_{\alpha^r}(L_{1}^{5}) $& Type of $ \varGamma_{\alpha^0}(L_{1}^{5}) $    &$ Der_{\alpha^r}(L_{1}^{5}) $&CN\\
		\hline
		$ r=0$&$ \begin{pmatrix}
		c_1&c_2\\
		0&c_1
		\end{pmatrix} $	&not  small   &$ \begin{pmatrix}
		0	&0\\
		0&0
		\end{pmatrix} $& Yes\\		
		\hline
		$ r\geq 1$&$ \begin{pmatrix}
0&c_2\\
0&0
\end{pmatrix} $	&  &$ \begin{pmatrix}
0	&d_2\\
0&0
\end{pmatrix} $& \\		
\hline		
	\end{tabular}
\end{center}
%%%%%%%%%%%%%%%%%%%%%
\begin{align*}
L_{3}^{6}:[e_{1},e_{1}]&=0,&     [e_{1},e_{2}]&=0,
&
[e_2,e_{1}]&=0
,&       [e_2,e_{2}]&=e_1,&\\
\alpha(e_{1})&=0,& \alpha(e_{2})&=e_{1}&&& &&   	
\end{align*} 	
\begin{center}
	\begin{tabular}{|c|c|c|c|c|}	
		\hline
		$ \alpha^r $ &$ \varGamma_{\alpha^r}(L_{3}^{6}) $& Type of $ \varGamma_{\alpha^0}(L_{3}^{6}) $    &$ Der_{\alpha^r}(L_{3}^{6}) $&CN\\
		\hline
		$ r=0$&$ \begin{pmatrix}
		c_1&c_2\\
		0&c_1
		\end{pmatrix} $	&  small   &$ \begin{pmatrix}
		0	&d_2\\
		0&0
		\end{pmatrix} $& Yes\\		
		\hline
		$ r\geq 1$&$ \begin{pmatrix}
		0&c_2\\
		0&0
		\end{pmatrix} $	&  &$ \begin{pmatrix}
		0	&d_2\\
		0&0
		\end{pmatrix} $& \\		
		\hline		
	\end{tabular}
\end{center}
%%%%%%%%%%%%%%%%%%%%%
\begin{align*}
L_{1}^{7}:[e_{1},e_{1}]&=0,&     [e_{1},e_{2}]&=0,
&
[e_2,e_{1}]&=0
,&       [e_2,e_{2}]&=e_1,&\\
\alpha(e_{1})&=e_{1},& \alpha(e_{2})&=e_{1}+e_{2}&&& &&   	
\end{align*} 	
\begin{center}
	\begin{tabular}{|c|c|c|c|c|}	
		\hline
		$ \alpha^r $ &$ \varGamma_{\alpha^r}(L_{1}^{7}) $& Type of $ \varGamma_{\alpha^0}(L_{1}^{7}) $    &$ Der_{\alpha^r}(L_{1}^{7}) $&CN\\
		\hline
		$ r=0$&$ \begin{pmatrix}
		c_1&c_2\\
		0&c_1
		\end{pmatrix} $	&  small   &$ \begin{pmatrix}
		0	&d_2\\
		0&0
		\end{pmatrix} $& Yes\\		
		\hline
		$ r\geq 1$&$ \begin{pmatrix}
		c_1&c_2\\
		0&c_1
		\end{pmatrix} $	&  &$ \begin{pmatrix}
		0	&d_2\\
		0&0
		\end{pmatrix} $& \\		
		\hline		
	\end{tabular}
\end{center}

%Il faut verifier les cas $ r=0 $. En cas de probleme remplacer par $ r\in \N^* $.\\
Summarizing the obtained results, the dimensions of the spaces of $\alpha^r$-derivations of Hom-Leibniz algebras in dimension $2$ vary between $0, 1$ and $2$.

%$ \dim( Der_{\alpha^r}(L_{i}^{j}))\leq 2. $\\
We have $ \dim(Der_{\alpha^r}(L_{i}^{j}))=2 $ in the following cases:
\begin{enumerate}[label=\upshape{\arabic*)},leftmargin=30pt]
\item $  Der_{\alpha^r}(L_{1}^{2})= Der_{\alpha^r}(L_{3}^{2})=<\begin{pmatrix}
1&0\\
0&0
\end{pmatrix},\begin{pmatrix}
0&0\\
0&1
\end{pmatrix}>  $, with $\alpha=\begin{pmatrix}
0&0\\
0&1
\end{pmatrix},$
  \item $  Der_{\alpha^r}(L_{2}^{1}) =<\begin{pmatrix}
2&0\\
0&1
\end{pmatrix},\begin{pmatrix}
0&1\\
0&0
\end{pmatrix}>  $, with $\alpha=id,$
  \item $  Der_{\alpha^r}(L_{4}^{1}) =<\begin{pmatrix}
1&0\\
0&0
\end{pmatrix},\begin{pmatrix}
0&1\\
0&0
\end{pmatrix}>  $, with $ \alpha =id$.
\end{enumerate}

%%%%%%%%%%%%%%%%%%%%%%%%%%%%%%%%%%
%\textbf{Voir degre de nilpotence?}\\
%\textbf{Voir simple ou non?}\\
%\textbf{Voir les autres articles?}\\

For $ \dim(Der_{\alpha^r}(L_{i}^{j}))=0 $, we distinguish two cases:

\begin{enumerate}[label=\upshape{\arabic*)},leftmargin=30pt]
  \item $ Der_{\alpha^r}(L_{1}^{5})$, $ r=0 $, $ z_1\neq -1 $,
%  ($ r=0 $?),
with $\alpha=\begin{pmatrix}
b^2&0\\
0&b
\end{pmatrix}$,
  \item $ Der_{\alpha^r}(L_{2}^{6}) $, $ r=0 $ with $\alpha=\begin{pmatrix}
0&1\\
0&0
\end{pmatrix}.$
\end{enumerate}

For the other cases of $Der_{\alpha^r}(L_{i}^{j})$, the dimension is equal to $1$.

Moreover, $ \dim(\varGamma_{\alpha^r}(L_{i}^{j}))$ vary between one and two.

We have $\dim(\varGamma_{\alpha^r}(L_{i}^{j})) =2$ for the following cases:
\begin{enumerate}[label=\upshape{\arabic*)},leftmargin=30pt]
  \item $\varGamma_{\alpha^r}(L_{1}^{1})=\varGamma_{\alpha^r}(L_{2}^{1})=\varGamma_{\alpha^r}(L_{3}^{1})=<id, \begin{pmatrix}
0&1\\
0&0
\end{pmatrix}> $, $ \alpha=id, $
  \item $\varGamma_{\alpha^r}(L_{2}^{6})=\varGamma_{\alpha^r}(L_{3}^{6})=<id, \begin{pmatrix}
0&1\\
0&0
\end{pmatrix}> $, $ \alpha=\begin{pmatrix}
0&1\\
0&0
\end{pmatrix},$
  \item $\varGamma_{\alpha^r}(L_{1}^{7})=<id, \begin{pmatrix}
0&1\\
0&0
\end{pmatrix}> $, $ \alpha=\begin{pmatrix}
1&1\\
0&1
\end{pmatrix}.$
\end{enumerate}

%The other case $ 1 $.
%%\begin{enumerate}[(i)]
%%	\item
%%	\item
%%	\item
%%\end{enumerate}
%%\textbf{Derivation dim1}
For the other cases of $\dim(\varGamma_{\alpha^r}(L_{i}^{j}))$ the dimension is equal to $1$.\\

%%%%%%%%%%%%%%%%%%%%%%%%%%%%%%%%%%%%%%%%%%%%%%%%%%%%%%%%%%%%%
%%%%%%%%%%%%%%%%%%%%%%%%%%%%%%%%%%%%%%%%%%%%%B3B3B3B3333333333333333333333333333333333333333333333
%%%%%%%%%%%%%%%%%%%%%%%%%%%%%%%%%%%%%%%%%%%%%%%%%%%%%%%%%%%%%%%%%%%%%%%%%%%%%%%%%%

\end{document}